\magnification=1200
\documentstyle{amsppt}
\document
\centerline{THE HOMOLOGY OF ITERATED LOOP SPACES}
\vskip 6pt
\centerline{\sl V.A.Smirnov\footnote"*"{Supported by RFFI Grant 
99-01-00114}}
\centerline{APPENDIX of  {\sl F. Sergeraert}}
\vskip .5cm
\centerline{Introduction}
\vskip 6pt
In recent years,  to solve various problems in Algebraic  Topology  it
has been necessary to consider more and more complicated structures on
the singular chain complex $C_*(X)$ of a topological space $X$ and its
homology $H_*(X)$.

One of  the  most  difficult problem is the problem of calculating the
homology groups of  iterated  loop  spaces.  The  first  steps  toward
solving  this  problem were made by J.F.Adams,  [1].  To calculate the
homology $H_*(\Omega X)$ of the loop space $\Omega X$ of a topological
space  $X$ he introduced the notion of the cobar construction $FK$ on
a coalgebra $K$.

Recall that a chain complex $K$ is called  a  coalgebra  if  there  is
given  a  chain  mapping $\nabla\colon K\to K\otimes K$ satisfying the
associativity                relation
$$(\nabla\otimes 1)\nabla=(1\otimes\nabla)\nabla.$$

The cobar  construction  $FK$  of  a  coalgebra  $K$ is a differential
algebra that coincides,  as a graded algebra,  with the tensor algebra
$TS^{-1}K$ on  the  desuspension  $S^{-1}K$ of $K$.  The generators in
$FK$ are denoted by $[x_1,\dots,x_n]$ where $x_i\in K$, $1\le i\le n$,
and have dimensions $$dim[x_1,\dots,x_n]=\sum_{i=1}^ndim(x_i)-n.$$

The product $\pi\colon FK\otimes FK\to FK$ is defined by the formula
$$\pi([x_1,\dots,x_n]\otimes [x_{n+1},\dots,x_{n+m}])=
[x_1,\dots,x_{n+m}].$$

The differential on the generators $[x]\in FK$ is defined by the
formula $$d[x]=-[d(x)]+\sum(-1)^\epsilon [x',x''],$$ where
$\nabla(x)=\sum x'\otimes x''$, $\epsilon =dim(x')$. On the other
elements, the differential is defined as a graded derivation.

The chain complex $C_*(X)$ of a topological space $X$ possesses a
natural coalgebra structure. Therefore there is defined the cobar
construction $FC_*(X)$.

J.F.Adams proved that for a  simply-connected  topological  space  $X$
there  is  a  chain  equivalence of differential algebras $$C_*(\Omega
X)\simeq FC_*(X).$$

In particular,  if the topological space $X$ is the  suspension  of  a
space $Y$,  i.e. $X=SY$, then the coalgebra structure on $C_*(X)$ is
can be taken as trivial.  Hence the chain complex $C_*(\Omega  X)$  of
the loop  space  $\Omega  X$  will  be  chain equivalent to the tensor
algebra $TC_*(Y)$ on the chain complex  $C_*(Y)$.

Unfortunately, the cobar construction $FK$ of a coalgebra $K$ does not
admit  iteration,  because passing to the cobar construction we lose a
coalgebra structure.  There  is,  in  general,  no  natural  coalgebra
structure on the Adams cobar construction $FK$ on a coalgebra $K$.

However H.J.Baues [2] introduced a coalgebra structure on the Adams
cobar construction $FC_*(X)$  on  the  chain  complex  $C_*(X)$  of  a
topological space $X$. This structure was determined using a family of
operations $$\nabla_{n,m}\colon  C_*(X)\to  C_*(X)^{\otimes  n}\otimes
C_*(X)^{\otimes m}$$
of dimensions $n+m-1$. Thus he defined the double cobar construction
$F^2C_*(X)$ and proved that for a $2$-connected topological space $X$
there is a chain equivalence $$C_*(\Omega^2X)\simeq F^2C_*(X).$$

But there is no known appropriate coalgebra structure on the double
cobar construction and therefore further iterations are not possible.

In [3] J.P.May introduced the notion of an operad and investigated
the structure on iterated loop spaces. This structure is used in [4] to
calculate the    homology    $H_*(\Omega^nS^nX)$.   The   homology
$H_*(\Omega^nS^nX)$ has also been investigated in [5], [6] and others.

In [7], [8], [9] the operad methods were transfered from the category
of topological spaces to the category of chain complexes. It was shown
that on the singular chain complex $C_*(X)$ of a topological space $X$
there is a natural $E_\infty$-coalgebra structure.

Using this  structure,  the  chain  complex  $C_*(\Omega^nX)$  of  the
$n$-fold loop space $\Omega^nX$ of an $n$-connected topological  space
$X$  was expressed in terms of the chain complex $C_*(X)$ of the space
$X$.

Our aim here is to construct a spectral sequence for the  homology  of
iterated loop spaces and produce some calculations.

To do  it,  we  use the $E_\infty$-coalgebra structure on the singular
chain complex $C_*(X)$ of a topological space $X$ (Theorem 1). Then we
express the  chain  complex $C_*(\Omega^nX)$ of the iterated loop space
$\Omega^nX$ in terms of the cobar construction  of  $C_*(X)$  (Theorem
2).

After that we consider the spectral sequence of the cobar construction
and calculate its first term with $\Bbb Z/2$-coefficients (Theorem 3),
with  $\Bbb  Z/p$-coefficients  (Theorem  3')  and  over  a  field  of
characteristic zero (Theorem 3'').  Also we give the expression of the
$E^1$-term  as  the free $n$-Poisson algebra generated by the homology
of $X$ (Theorems 4, 4').

In the cases  of  $\Bbb  Z/2$,  $\Bbb  Z/p$  and  characteristic  zero
coefficients, we calculate the differential $d^1$ on the first term and
obtain the expression  of  the  $E^2$-term  of  the  spectral  sequence
(Theorems 5, 5', 5'').

Finally we  apply  these  results  to  calculate  the  homology of the
iterated loop spaces of the stunted real and complex projective spaces
(Theorem 6, 7, 8), which play important roles in Algebraic Topology.
Some calculations with these spaces have been produced by F.Cohen and
R. Levi [10].

Note that a general method for calculating the homology of iterated
loop spaces is  given by simplicial theory. There is the simplicial
construction $GX$ of the loop space of   a simplicial set $X$ and its
iteration $G^nX$.  But  this  construction  is  very  complicated  and
direct calculations may be produced only using computer methods [11].

The computer calculations may be produced for any topological space, but
only in low dimensions. Conversely the spectral sequence calculations
may be produced for ``nice" spaces in any dimensions.  So the computer
 calculations complement the spectral sequence calculations.

The problem of comparing the computer and the spectral sequence
calculations for these spaces was stated by F.Sergeraert.  It seems to
be very useful for both sides of these calculations.

\vskip .5cm
\centerline{\S 1. Operads and algebras over operads}
\vskip 6pt
Consider the category of chain complexes over a ring $R$. By a
symmetric family $\Cal E$ in this category is meant a family $\Cal
E=\{\Cal E(j)\}_{j\ge 1}$ of chain complexes $\Cal E(j)$ operated on
by the symmetric groups $\Sigma_j$.

Given two symmetric families $\Cal E$, $\Cal E'$ we define the
symmetric family   $\Cal   E\otimes\Cal   E'$   by   putting   $$(\Cal
E\otimes\Cal E')(j)=\Cal E(j)\otimes\Cal E'(j)$$ and
$(\Cal E\times\Cal E')$ by putting $(\Cal E\times\Cal E')(j)$  equal
to the quotient module of the $\Sigma_j$-free module generated by the
module
$$\sum_{j_1+\dots+j_k=j}\Cal E(k)\otimes\Cal E'(j_1)\otimes\dots
\otimes\Cal E'(j_k)$$
modulo the equivalence generated by the relations
$$\gather
x\sigma\otimes x'_1\otimes \dots\otimes x'_k\sim x'\otimes
x'_{\sigma^{-1}(1)}\otimes\dots\otimes x'_{\sigma^{-1}(k)}\cdot
\sigma(j_1,\dots,j_k),\\
x\otimes x'_1\sigma_1\otimes\dots\otimes x'_k\sigma_k\sim
x\otimes x'_1\otimes\dots\otimes x'_k\cdot\sigma_1\times\dots
\times\sigma_k.\endgather $$

Here $\sigma(j_1,\dots,j_k)$ is the permutation of a set of $j$
elements obtained by partitioning the set into $k$ blocks of
$j_1,\dots,j_k$ elements, respectively, and carrying out on these
blocks the permutation $\sigma$, while $\sigma_1\times\dots\times
\sigma_k$ means the image of the element $(\sigma_1,\dots,\sigma_k)$
under the imbedding $\Sigma_{j_1}\times\dots\times\Sigma_{j_k}\to
\Sigma_j$.

It is easy to see that for symmetric families $\Cal E,\Cal E',
\Cal F,\Cal F',$ there is the interchange mapping
$$T\colon (\Cal E\otimes\Cal E')\times(\Cal F\otimes\Cal E')\to
(\Cal E\times\Cal F)\otimes(\Cal E'\times\Cal F').$$

A symmetric family $\Cal E $ is called an {\it operad} if there is given
a symmetric-family mapping $\gamma\colon\Cal E\times\Cal E\to
\Cal E$ such that $\gamma(\gamma\times 1)=\gamma(1\times\gamma)$ or,
what is the same, the following diagram is commutative
$$\CD
\Cal E\times\Cal E\times\Cal E@>\gamma\times 1>>\Cal E\times\Cal E\\
@V1\times\gamma VV @VV\gamma V\\
\Cal E\times\Cal E@>\gamma >>\Cal E\endCD $$

If there exist an element $1\in\Cal E(1)$ such that $\gamma(1\otimes
x)=\gamma(x\otimes 1^{\otimes k})=x$ for all $x\in\Cal E(k)$ then
we say that $\Cal E$ is an operad with identity.

A mapping  of  operads  $f\colon\Cal  E\to\Cal  E'$  is a mapping of
symmetric families for which the following diagram is commutative
$$\CD \Cal E\times\Cal E@>\gamma>>\Cal E\\
@Vf\times fVV @VVfV\\
\Cal E'\times\Cal E'@>\gamma'>>\Cal E'\endCD $$

If $\Cal E$ and $\Cal E'$ are operads with identities $1$ and $1'$
respectively, then it is required that $f(1)=1'$.

We shall say that an operad $\Cal E$ acts on a symmetric family $\Cal
F$ on the left (right) if there is given a mapping $\mu'\colon
\Cal E\times\Cal F\to\Cal F$ ($\mu''\colon\Cal F\times\Cal E\to\Cal
F$) such that $$\mu'(\gamma\times 1)=\mu'(1\times\mu')\quad
(\mu''(1\times\gamma)=\mu''(\mu''\times 1))$$
or, what is the same, the following diagrams are commutative
$$\CD
\Cal E\times\Cal E\times\Cal F@>\gamma\times 1>>\Cal E\times\Cal F
\\@V1\times\mu'VV @VV\mu'V\\
\Cal E\times\Cal F@>\mu'>>\Cal F\endCD \qquad\left(\CD
\Cal F\times\Cal E\times\Cal E@>1\times\gamma >>\Cal F\times\Cal E\\
@V\mu''\times 1VV @VV\mu''V\\
\Cal F\times\Cal E@>\mu''>>\Cal F\endCD\right) $$

For any chain complex $X$ and symmetric family $\Cal E,$ we define
chain complexes $\Cal E(X)$, $\overline{\Cal E}(X)$ by putting
$$\Cal E(X)=\sum_k\Cal E(k)\otimes_{\Sigma_k}X^{\otimes k},
\quad\overline{\Cal E}(X)=\prod_kHom_{\Sigma_k}(\Cal E(k);X^{\otimes
k}).$$

If $\Cal E$ is an operad, then the operad structure in $\Cal E$
determines a mapping
$$\gamma\colon\Cal E^2(X)=\Cal E(\Cal E(X))\to\Cal E(X)$$
such that the correspondence $X\longmapsto \Cal E(X)$ is a monad
(also known as a triple) in the category of chain complexes.

Dually, the operad structure in $\Cal E$ determines a mapping
$$\overline\gamma\colon\overline{\Cal E}(X)\to\overline{\Cal E}^2
(X)=\overline{\Cal E}(\overline{\Cal E}(X))$$
such that the correspondence $X\longmapsto \overline{\Cal E}(X)$
is a comonad in the category of chain complexes.

A chain complex $X$ is called an {\it algebra over the operad} $\Cal E$,
or simply an $\Cal E$-agebra, if there is given a mapping $\mu\colon
\Cal E(X)\to X$ satisfying the associativity relation:
$$\mu\circ\gamma(X)=\mu\circ\Cal E(\mu)$$ or, what is the same,
the following diagram is commutative
$$\CD \Cal E^2(X)@>\gamma(X)>>\Cal E(X)\\
@V\Cal E(\mu)VV @VV\mu V\\
\Cal E(X)@>\mu >>X\endCD $$

Dually, a chain complex $X$ is called a {\it coalgebra over the operad}
$\Cal E$, or simply an $\Cal E$-coalgebra, if there is given a mapping
$\tau\colon X\to\overline{\Cal E}(X)$ satisfying the associativity
relation: $$\overline\gamma(X)\circ\tau=\overline{\Cal E}(\tau)
\circ\tau$$ or, what is the same, the following diagram is commutative
$$\CD X@>\tau >>\overline{\Cal E}(X)\\
@V\tau VV @VV\overline\gamma(X)V\\
\overline{\Cal E}(X)@>\overline{\Cal E}(\tau)>>\overline{\Cal E}^2(X)
\endCD $$

Consider some examples of operads and algebras over operads.

1. The  operad $E_0=\{E_0(j)\}$,  where $E_0(j)$ -- the free
$R$-module with one zero  dimensional  generator  $e(j)$  and  trivial
action  of  the  symmetric group $\Sigma_j$.  So $E_0(j)\cong R$.  The
operation $\gamma\colon E_0\times E_0\to E_0$ is given by the  formula
$$\gamma(e(k)\otimes        e(j_1)\otimes\dots\otimes       e(j_k))=
e(j_1+\dots+j_k).$$

It is easy to see that the required relations are satisfied and
algebras (coalgebras) over the operad $E_0$ are simply
commutative and associative algebras (coalgebras).

2. The operad $A=\{A(j)\}$, where $A(j)$ is the $\Sigma_j$-free
module
with one zero dimensional generator $a(j)$. So $A(j)\cong R(\Sigma_j)$.
The operation $\gamma\colon A\times A\to A$ is given by the formula
$$\gamma(a(k)\otimes a(j_1)\otimes\dots\otimes a(j_k))=
a(j_1+\dots+j_k).$$

It is easy to see that the required relations are satisfied and
that algebras (coalgebras) over the operad $A$ are simply
associative algebras (coalgebras).

3. For a symmetric family $\Cal E,$ define the suspension $S\Cal E$
by putting $(S\Cal E)(j)=S^{j-1}\Cal E(j),$ the $(j-1)$-fold
suspension over $\Cal E(j)$. It is clear that if $\Cal E$ is an operad
then the suspension $S\Cal E$ will also be an operad, and if $X$ is an
algebra (coalgebra) over an operad $\Cal E$ then the  suspension  $SX$
will also be an algebra (coalgebra) over the operad $S\Cal E$.

4. For operads $\Cal E$, $\Cal E',$ the tensor product $\Cal E \otimes
\Cal E'$ evidently is an operad.

5. Let $L$ be the suboperad of the operad $A$ generated by the element
$$b(2)=a(2)-a(2)T,\quad T\in\Sigma_2.$$ Then algebras (coalgebras)
over the operad $L$ will be simply Lie algebras (coalgebras).

Similary, let $L_n$ be the suboperad of the operad $S^nA$ generated by
the element $$b_n(2)=s^na(2)+(-1)^ns^na(2)T.$$ Then algebras
(coalgebras) over the operad $L_n$ will be simply $n$-Lie algebras
(coalgebras) [4]. It means there is given a Lie bracket of
dimension $n$, called an $n$-Lie bracket satisfying the relations
$$\gather
[x,y]+(-1)^\epsilon [y,x]=0;\\
[x,[y,z]]=[[x,y],z]+(-1)^\epsilon [y,[x,z]];\endgather $$
where $\epsilon =(dim(x)+n)\cdot (dim(y)+n)$.
So this is defferent from $N$-Lie meaning $N$-ary bracket [12].

6. Let $P_n=E_0\times L_n$ and the operad structure $\gamma$ is
determined by the corresponding structures in $E_0$, $L_n$ and
by the formulas
$$\gamma(b_n(2)\otimes 1\otimes e(2))=e(2)\otimes b_n(2)\otimes 1
+e(2)\otimes 1\otimes b_n(2)\cdot (213).$$
The operad   $P_n$   is   called   the  $n$-Poisson  operad.  Algebras
(coalgebras)   over   this   operad   called   $n$-Poisson    algebras
(coalgebras). They  are  commutative  algebras  together  with  a  Lie
bracket of dimension $n$ satisfying the Poisson relation
$$[x,y\cdot z]=[x,y]z+(-1)^\delta y[x,z],\quad [4],$$
where $\delta=(dim(x)+n)\cdot dim(y)$.

7. For any chain complex $X$ define the  operads  $\Cal  E_X$,  $\Cal
E^X$,  by  putting  $$\Cal  E_X(j)=Hom(X^{\otimes j};X);\quad \Cal
E^X(j)=Hom(X;X^{\otimes j}).$$

The actions of the symmetric groups are determined by the permutations
of factors of $X^{\otimes j}$ and the operad structures are defined by
the formulas
$$\gather \gamma_X(f\otimes g_1\otimes\dots\otimes g_k)=
f\circ(g_1\otimes\dots\otimes g_k),\quad f\in\Cal E_X(k),~
g_i\in\Cal E_X(j_i);\\
\gamma^X(f\otimes g_1\otimes\dots\otimes g_k)=
(g_1\otimes\dots\otimes g_k)\circ f,\quad f\in \Cal E^X(k),~
g_i\in\Cal E^X(j_i).\endgather$$

Directly from the definitions, it follows that a chain complex $X$
is an algebra (coalgebra) over the operad $\Cal E$ if and only if
there is given an operad mapping $\xi\colon\Cal E\to\Cal E_X$
($\xi\colon \Cal E\to\Cal E^X$).

Analogously, for chain complexes $X$ and $Y$ there are defined the
symmetric families $\Cal F_{X,Y}$, $\Cal F^{X,Y}$
$$\Cal F_{X,Y}(j)=Hom(X^{\otimes j};Y),\quad
\Cal F^{X,Y}(j)=Hom(X;Y^{\otimes j})$$ and actions
$$\gather \mu'\colon\Cal E_Y\times\Cal F_{X,Y}\to\Cal F_{X,Y},~
\mu''\colon\Cal F_{X,Y}\times\Cal E_X\to\Cal F_{X,Y};\\
\mu'\colon\Cal E^X\times\Cal F^{X,Y}\to\Cal F^{X,Y},~
\mu''\colon\Cal F^{X,Y}\times\Cal E^Y\to\Cal F^{X,Y}.\endgather $$

8. One of the most important topological operad is the
J.M.Boardman and R.M.Vogt's ``little $n$-cubes'' operad $E_n$, [13].

Let $J$ denote the
open interval $(0,1)$ and $J^n$ the open $n$-dimensional cube.
By an $n$-dimensional little cube is meant  an affine embedding
$f\colon J^n\to J^n$ with parallel axes. Then $E_n(j)$ is the set of
ordered $j$-tuples $(f^1,\dots,f^j)$ of $n$-dimensional little cubes
$f^i\colon J^n\to J^n$ such that images don't intersect.
This operad acts on the $n$-fold loop space $\Omega^nX$ of   a space
$X$.

The direct limit of the operads $E_n$ over the inclusions $E_n
\subset E_{n+1}$ is denoted $E_\infty$. It is an acyclic operad
with free actions of the symmetric groups.

9. It is easy to see that if $\Cal E$ is a topological operad,
then its singular chain complex $C_*(\Cal E)$ is an operad in
the category of chain complexes, and if $\Cal E$ acts on a space
$X$ then $C_*(\Cal E)$  acts on $C_*(X)$.

Similary, the homology $H_*(\Cal E)$ of a topological operad $\Cal E$
is an operad in the category of graded modules, and if $\Cal E$ acts
on a space $X,$ then $H_*(\Cal E)$ acts on the homology $H_*(X)$.

In particular, the homology $H_*(E_n)$ of the topological $n$-cubes
operad $E_n$ is isomorphic to the $n$-Poisson operad $P_n$, [15].

So the homology $H_*(\Omega^nX)$ of the $n$-fold loop space
$\Omega^nX$, $n>1$ is an algebra over the $n$-Poisson operad $P_n$.

The operad $C_*(E_\infty)$ gives us an example of an acyclic operad
in the category of chain complexes with free actions of the symmetric
groups.

Note that all acyclic operads with the free actions of the symmetric
groups  consist of $\Sigma_j$-homotopy equivalent chain complexes.  We
will call such operads: $E_\infty$-operads.

10. Another example of an $E_\infty $-operad  is  given  simplicial
resolutions of the symmetric groups. Denote by $E\Sigma_*(j)$ the free
simplicial resolution of the symmetric group $\Sigma_j$, i.e.
$$E\Sigma_*(j):\Sigma_j@<<<\Sigma_j\times\Sigma_j@<<<\dots $$
The mappings $\gamma\colon\Sigma_k\times\Sigma_{j_1}\times\dots
\times\Sigma_{j_k}\to\Sigma_{j_1+\dots+j_k}$ induce the operad
structure
$$\gamma_*\colon E\Sigma_*(k)\times E\Sigma_*(j_1)\times\dots\times
E\Sigma_*(j_k)\to E\Sigma_*(j_1+\dots+j_k).$$
So $E\Sigma_*$ will be an acyclic operad with free actions of the
symmetric groups in the category of simplicial sets.

Taking the  chain   complex   $C_*(E\Sigma_*),$   we   obtain   an
$E_\infty$-operad in the category of chain complexes. Denote it simply
by $E\Sigma$.

Note that for any chain operad $\Cal E,$ the operad $\Cal E\otimes
E\Sigma$ has the same homology and free actions of the
symmetric groups. The projection $E\Sigma\to E_0$ induces the
projection $\Cal E\otimes E\Sigma\to \Cal E$. So $\Cal E\otimes
E\Sigma$ may be considered as a  $\Sigma$-free resolution of the operad
$\Cal E$. If $\Cal E$ is an acyclic operad, then $\Cal E\otimes E\Sigma$
will be an  $E_\infty$-operad.

\vskip .5cm
\centerline{\S 2. On the chain complex of a topological space}
\vskip 6pt
Here we consider structure on the singular chain complex $C_*(X)$
of a topological space $X$, and dually on the singular cochain complex
$C^*(X)$.

Besides the coalgebra structure $$\nabla\colon C_*(X)\to C_*(X)
\otimes C_*(X)$$ on the chain complex of a topological space, there are
coproducts $$\nabla_i\colon C_*(X)\to C_*(X)\otimes C_*(X)$$
increasing dimensions by $i$ and such that
$$d(\nabla_i)=\nabla_{i-1} +(-1)^iT\nabla_{i-1},$$   where  $T\colon
C_*(X)\otimes C_*(X)\to C_*(X)\otimes C_*(X)$ permutes factors.

Dually, on the cochain complex $C^*(X)$ besides the algebra structure
$$\cup\colon C^*(X)\otimes C^*(X)\to C^*(X)$$ there are products
$$\cup_i\colon C^*(X)\otimes C^*(X)\to C^*(X)$$ such that
$d(\cup_i)=\cup_{i-1}+(-1)^i\cup_{i-1}T$.

To describe all operations on the singular chain complex $C_*(X)$ of a
toplogical space $X$ and on its dual cochain complex $C^*(X),$ we
consider the corresponding operad.

For $n\ge 0$ denote by $\Delta^n$ the normalized chain complex of the
standard $n$-dimensional simplex. Then $\Delta^*=\{\Delta^n\}$ is a
cosimplicial object in the category of chain complexes. Consider also
the cosimplicial object $(\Delta^*)^{\otimes j}=\Delta^*\otimes\dots
\otimes\Delta^*$ and
$$E^\Delta(j)=Hom(\Delta^*;(\Delta^*)^{\otimes j}),$$
where $Hom$ is considered in the category of cosimplicial objects.

The family $E^\Delta=\{E^\Delta(j)\}$ will be the operad for which
the actions of the symmetric groups and the operad structure are
defined similary to the corresponding structures for the above defined
operad $\Cal E^X$, where instead of $X$ we take $\Delta^*$.

Note that the complexes $\Delta^n$ are acyclic and hence the operad
$E^\Delta$ is also acyclic.

{\bf Theorem 1.} {\sl On the chain complex $C_*(X)$ of a topological
space $X$ there exists a natural  $E^\Delta$-coalgebra structure
$\tau\colon
C_*(X)\to\overline E^\Delta(C_*(X))$ with the following universal
property: if for some operad $\Cal E$ there is a natural
$\Cal E$-coalgebra structure $\widetilde\tau\colon C_*(X)\to
\overline{\Cal E}(C_*(X))$, then there exist a unique operad
mapping $\xi\colon\Cal E\to E^\Delta$ such that the following
diagram commutes
$$\CD C_*(X)@>\tau>>\overline E^\Delta(C_*(X))\\
@V=VV @VV\overline\xi V\\
C_*(X)@>\widetilde\tau>>\overline{\Cal E}(C_*(X))\endCD $$}

{\sl Proof.} Our aim is to define natural operations
$$\tau\colon E^\Delta(j)\otimes C_*(X)\to C_*(X)^{\otimes j}.$$
Let $x_n\in C_n(X)$, $y\in E^\Delta(j)=Hom(\Delta^*;(\Delta^*)
^{\otimes j})$. The element $x_n\in C_n(X)$ determines a chain
mapping $\overline x_n\colon\Delta^n\to C_*(X)$ such that the generator
$u_n\in\Delta^n$ maps to $x_n$.

Note that there is the operation $\tau^n\colon E^\Delta(j)\otimes
\Delta^n\to(\Delta^n)^{\otimes j}$. Define the required operation $\tau$
by putting $$\tau(y\otimes x_n)=(\overline x_n)^{\otimes j}\circ\tau^n
(y\otimes u_n).$$
Then we have the following commutative diagram
$$\CD E^\Delta(j)\otimes C_*(X)@>\tau >> C_*(X)^{\otimes j}\\
@A1\otimes \overline x_nAA @AA(\overline x_n)^{\otimes j}A\\
E^\Delta(j)\otimes\Delta^n @>\tau^n>>(\Delta^n)^{\otimes j}\endCD$$
It is easy to see that the required relations are satisfied.

An  $E^\Delta$-coalgebra structure on the chain complex $C_*(X)$
of a topological space $X$ induces  an  $E^\Delta$-algebra structure
on the cochain complex $$C^*(X)=Hom(C_*(X);R).$$ The corresponding
operations $\mu\colon E^\Delta(j)\otimes C^*(X)^{\otimes j}\to C^*(X)$
are defined by the formulas
$$\mu(y\otimes f_1\otimes\dots\otimes f_j)(x)=(f_1\otimes\dots\otimes
f_j)\circ\tau(y\otimes x),$$
where $y\in E^\Delta(j)$, $f_i\colon C_*(X)\to R$, $x\in C_*(X)$.
So we have

{\bf Theorem 1'.} {\sl On the cochain complex $C^*(X)$ of a topological
space $X$ there exists a natural $E^\Delta$-algebra structure $\mu\colon
E^\Delta(C^*(X))\to C^*(X)$ with the following universal property:
if for some operad $\Cal E,$ there is a natural $\Cal E$-algebra
structure $\widetilde\mu\colon\Cal E(C^*(X))\to C^*(X)$, then there
exists a unique operad mapping $\xi\colon\Cal E\to E^\Delta$ such
that the following diagram commutes
$$\CD \Cal E(C^*(X))@>\widetilde\mu>>C^*(X)\\
@V\xi VV @VV=V\\
E^\Delta(C^*(X))@>\mu>>C^*(X)\endCD $$}

Let $R(\Sigma_2)$ be the $\Sigma_2$-free resolution with generators
$e_i$ of dimensions $i$ and differential defined by the formula
$$d(e_i)=e_{i-1}+(-1)^ie_{i-1}T,\quad T\in\Sigma_2.$$

Since $E^\Delta(2)$ is acyclic, there is a $\Sigma_2$-chain mapping
$R(\Sigma_2)\to E^\Delta(2)$ and hence a mapping $$R(\Sigma_2)
\otimes_{\Sigma_2}C^*(X)^{\otimes 2}\to C^*(X).$$ Its restriction on
the elements $e_i$ is usually denoted by $$\cup_i\colon C^*(X)\otimes
C^*(X)\to C^*(X)$$ and called cup-$i$ product.

Let $p^n\colon\Delta^n\to S\Delta^{n-1}$ be the projection obtained by
contracting  the ($n-1$)-dimensional face spanned by the vertices with
numbers $0,1,\dots,n-1$.  These projections induce the  projection  of
operads  $E^\Delta\to  SE^\Delta$.  For  a  topological space $X,$ the
suspension $SC_*(X)$ will be a coalgebra over the operad  $SE^\Delta.$
The          following          diagram         commutes
$$\CD
SC_*(X)@>>>\overline{SE}^\Delta(SC_*(X))\\
@VVV @VVV\\
C_*(SX)@>>>\overline E^\Delta(C_*(SX)).\endCD $$

Iterating this construction we obtain the projections
$E^\Delta\to S^nE^\Delta$ and the commutative diagrams
$$\CD S^nC_*(X)@>>>\overline{S^nE}^\Delta(S^nC_*(X))\\
@VVV @VVV\\
C_*(S^nX)@>>>\overline E^\Delta(C_*(S^nX)).\endCD $$

Note that the operad $E^\Delta$ may be not $\Sigma$-free and so it
is not an  $E_\infty$-operad. To obtain an  $E_\infty$-operad we
consider the operad $E^\Delta\otimes E\Sigma$ and denote it simply by
$E$.

The projection $E\to E^\Delta$ induces an  $E$-coalgebra structure
on the chain complex $C_*(X)$ of a topological space $X$. The mapping
$E^\Delta\to SE^\Delta$ induces the operad mapping  $E\to SE$ and
for the chain complex $C_*(S^nX)$ there is the corresponding commutative
diagram similar to the diagram for the operad $E^\Delta$.

\vskip .5cm
\centerline{\S 3. Bar and cobar constructions over operads}
\vskip 6pt
Let $\Cal E$ be an operad with right action $\nu\colon\Cal F\times\Cal
E\to\Cal F$ on a symmetric family $\Cal F$, $\nu\colon\Cal F\times\Cal
E\to\Cal F$, and let $X$ be an algebra over the operad $\Cal E$, given
by $\mu\colon\Cal E(X)\to X$.

Consider the simplicial object (monadic  bar  construction)  $$B_*(\Cal
F,\Cal  E,X)=\{B_n(\Cal  F,\Cal E,X)\}$$ for which $B_n(\Cal F,\Cal
E,X)=\Cal F\Cal E^n(X)$ with face and degeneracy operators given by
the formulas
$$\gather d_0=\nu 1^n,~d_i=1^i\gamma 1^{n-i},~0<i<n;\\
d_n=1^n\mu,~s_j=1^{j+1}i1^{n-j+1},0\le j\le n.
\endgather $$

Its realization
$$B(\Cal F,\Cal E,X)=|B_*(\Cal F,\Cal E,X)|=\sum_n\Delta^n
\otimes B_n(\Cal F,\Cal E,X)/\sim $$
is called the bar construction.

In the case of trivial $\Cal F$, i.e. $\Cal F(1)=R$ and $\Cal F(j)
=0$ if $j\ge 2$, the corresponding bar construction is denoted by
$B(\Cal E,X)$.

If $X$ is an $A$-algebra, i.e. simply an algebra, then the bar
construction $B(A,X)$ will be chain equivalent to the desuspension  of
the  usual  bar  construction $B(X)$,  i.e.  $$B(A,X)\simeq S^{-1}BX,~
[8].$$

Dually, if $X$ is a coalgebra $\tau\colon X\to^M \overline{\Cal E}(X)$
over the operad $\Cal E$, then we can consider the cosimplicial object
$F^*(\Cal F,\Cal E,X)=\{F^n(\Cal F,\Cal E,X)\}$ for which  $F^n(\Cal
F,\Cal  E,X)=\overline{\Cal  F}\overline{\Cal E}^n(X)$ with  coface
and codegeneracy operators given by the formulas
$$\gather \delta^0=\overline\nu 1^n,~\delta^i=1^i\overline\gamma
1^{n-i},~0<i<n;\\
\delta^n=1^n\tau,~\sigma^j=1^{j+1}p1^{n-j+1},0\le j\le n.\endgather $$

Its realization
$$F(\Cal F,\Cal E,X)=|F^*(\Cal F,\Cal E,X)|=Hom(\Delta^*;F^*(\Cal
F,\Cal E,X)),$$ where $Hom$ is considered in the category of cosimplicial
objects, is called the cobar construction.

In the case of trivial $\Cal F$ the corresponding bar construction
is denoted by $F(\Cal E,X)$.

If $X$ is an $A$-coalgebra, i.e. simply a coalgebra then the cobar
construction $F(A,X)$ will be chain equivalent to the suspension of
the usual Adams cobar construction $F(X)$, i.e. $$F(A,X)\simeq SFX,~ [8].$$

Let now $X$ be an $n$-connected topological space. As was shown above,
on the chain complex $S^nC_*(\Omega^nX)$ there is the $S^nE$-coalgebra
structure
$$\tau\colon S^nC_*(\Omega^nX)\to\overline{S^nE}(S^nC_*(\Omega^nX)).$$

This structure  and the mapping $j\colon S^n\Omega^nX\to X$ induce the
mapping $$S^nC_*(\Omega^nX)\to\overline{S^nE}(S^nC_*(\Omega^nX))\to
\overline{S^nE}(C_*(X)).$$ This mapping is a coaugmentation of the
cosimplicial object $F^*(S^nE,E,C_*(X))$ and hence it induces a
mapping of $S^nE$-coalgebras: $$S^nC_*(\Omega^nX)\to F(S^nE,E,C_*(X)).$$

{\bf Theorem 2.} {\sl For any $n$-connected topological space $X,$ the
mapping
$$S^nC_*(\Omega^nX)\to F(S^nE,E,C_*(X))$$ is a chain equivalence of
$S^nE$-coalgebras.}

{\sl Proof.} For $n=1$ the chain equivalence
$SC_*(\Omega X)\to F(SE,E,C_*(X))$ follows from the Adams chain
equivalence $$SC_*(\Omega X)\simeq SFC_*(X)$$
and from the chain equivalence $$SFC_*(X)\simeq F(SE,E,C_*(X)),~[8].$$
Suppose that for any ($n-1$)-connected topological space $X,$ we have
a chain equivalence  $$S^{n-1}C_*(\Omega^{n-1}X)\to
F(S^{n-1}E,E,C_*(X)).$$

Then for an $n$-connected topological space $X$ we will have the
following sequence of chain equivalences
$$\gather S^nC_*(\Omega^nX)\simeq SF(S^{n-1}E,E,C_*(\Omega X))
\simeq F(S^nE,SE,SC_*(\Omega X))\simeq \\
F(S^nE,SE,F(SE,E,C_*(X)))\simeq F(S^nE,E,C_*(X)).\endgather$$
The composition will be the desirable chain equivalence
$$S^nC_*(\Omega^nX)\simeq F(S^nE,E,C_*(X)).$$

\vskip .5cm
\centerline{\S 4. A   spectral sequence for the homology of iterated
loop spaces}
\vskip 6pt
Let $X$ be an $E$-coalgebra. Consider the spectral sequence of the
cobar construction $F(S^nE,E,X)$ with respect to the filtration
determined by the operad grading.  Namely,  for any operad $\Cal E$ we
define the grading of mappings $f\colon \Cal E(j)\to X^{\otimes j}$ to
be equal  to $j$.

Thus we will have the grading of the elements of $\overline{S^nE}(X)$,
$\overline{S^nE}\circ\overline E(X)$ and so on. Similarly we will have
the grading of the elements of the  cobar  construction  $F(S^nE,E,X)$
which determines in it a decreasing filtration.

Note that the mappings
$\overline\mu\colon\overline{S^nE}\to\overline{S^nE}\circ\overline E$,
$\overline\gamma\colon\overline E\to\overline E\circ\overline E$
preserve gradings and the mapping $\tau\colon X\to \overline E(X)$
increases it. Therefore the first term of the corresponding spectral
sequence will be isomorphic  to the homology of the cobar construction
$F(S^nE,E,X)$, where $X$ is considered as a trivial $E$-coalgebra.

To calculate its homology, recall the notion of the iterated
cobar construction over a cocommutative coalgebra $K$.

The cobar construction $FK$ will be a cocommutative Hopf algebra  with
a coproduct   $\nabla\colon   FK\to   FK\otimes  FK$  defined  on  the
generators $[x]\in FK$ by the formula
$$\nabla[x]=[x]\otimes 1+1\otimes [x].$$

Note that  a  commutative  coalgebra  $K$  may  be  considered as an
$E$-coalgebra.  The required $E$-coalgebra structure is induced by the
projection $E\to  E_0$.
Moreover, the chain equivalence $SFK\simeq F(SE,E,K)$ will be a chain
equivalence of $SE$-coalgebras.

Thus, in this case the Adams cobar construction may be iterated and
there are chain equivalences $$S^nF^nK\simeq F(S^nE,E,K).$$
In particular, if a chain complex $X$ has a   trivial $E$-coalgebra
structure then there are the chain equivalences
$$S^nF^nX\simeq F(S^nE,E,X).$$
So to calculate the homology of the cobar construction $F(S^nE,E,X)$
of   a trivial $E$-coalgebra $X,$ it is sufficient to calculate the
homology of the iterated Adams cobar construction $F^nX$.

Consider the case of $Z/2$-coefficients. The homology
$H_*(FX)$ is isomorphic to the tensor algebra $TS^{-1}H_*(X)$ on   the
desuspension $S^{-1}H_*(X)$ of   $H_*(X)$, i.e.
$$H_*(FX)\cong TS^{-1}H_*(X).$$

The double cobar construction $F^2X$ will be chain equivalent to the
cobar construction $FTS^{-1}H_*(X)$, where $TS^{-1}H_*(X)$ is considered
as a Hopf algebra with a coproduct
$$\nabla\colon TS^{-1}H_*(X)\to TS^{-1}H_*(X)\otimes TS^{-1}H_*(X)$$
determined on the generators $[x]\in S^{-1}H_*(X)$ by the formula
$$\nabla[x]=[x]\otimes 1+1\otimes [x].$$

For a  graded  module  $M$  (over  $Z/2$)  denote  by $LM$ the free Lie
algebra generated by $M$ and by $T_sM$  the  quotient  of  the  tensor
algebra $TM$ over the permutations of factors. Note that $T_sM$ is the
free commutative algebra generated by $M$.

The Hopf algebra structure on $TM$ induces the Hopf algebra  structure
on $T_sM$. By the Poincare-Birkhoff-Witt theorem there is an
isomorphism of coalgebras $$TM\cong  T_sLM$$  and  hence  there  is  an
isomorphism $$TS^{-1}H_*(X)\cong T_sLS^{-1}H_*(X).$$

For a  coalgebra  $K$  denote  by  $PK$  the  module of it's primitive
elements, i.e. $$PK=\{x\in K|\nabla(x)=x\otimes 1+1\otimes x\}.$$
Note that the module $PT_sM$ of primitive elements of the Hopf algebra
$T_sM$ is generated by the elements of the form $x^{2^k}$, $x\in M$,
$k\ge 0$.

The homology of $FT_sM$ is isomorphic to the free commutative
algebra generated by the module $S^{-1}PT_sM$ and hence  there  is  an
isomorphism
$$H_*(F^2X)\cong T_sS^{-1}PT_sLS^{-1}H_*(X).$$

So this construction may be iterated and by induction we obtain
isomorphisms
$$H_*(F^nX)\cong T_s(S^{-1}PT_s)^{n-1}LS^{-1}H_*(X),$$ and hence
isomorphisms
$$S^{-n}H_*(F(S^nE,E,X))\cong T_s(S^{-1}PT_s)^{n-1}LS^{-1}H_*(X).$$

Recall
that a graded module $L$ (over $Z/2$) is called an $n$-Lie algebra
if there is given an operation $[~,~]\colon L\otimes L\to L$ called
an $n$-Lie bracket of dimension $n$ and satisfying the relations
$$\gather
[x,x]=0;\\
[x,y]+[y,x]=0;\\
[x,[y,z]]=[[x,y],z]+[y,[x,z]].\endgather $$

For a graded module $M$ denote by $\Cal E_nM$
the module generated
by the elements $e_{i_1}\dots e_{i_k}x$, where $x\in M$, $0\le i_1\le
\dots \le i_k\le n$, and the dimensions of these elements are
defined equal to $i_1+2i_2+\dots+2^{k-1}i_k+2^kdim(x)$.

These elements $e_{i_1}\dots e_{i_k}x$ of $\Cal E_nM$ may be
rewritten using the Dyer-Lashof algebra $\Cal R$ in the form
$$Q^{j_1}\dots Q^{j_k}x;\quad j_1\le 2j_2,\dots,j_{k-1}\le 2j_k,
dim(x)\le j_k\le dim(x)+n,$$ where
$$\gather j_k=i_k+dim(x),\\j_{k-1}=i_{k-1}+i_k+2dim(x),\\
\dots \\
j_1=i_1+i_2+2i_3+\dots+2^{k-2}i_k+2^{k-1}dim(x).\endgather $$
It is clear that the sequences $Q^{j_1}\dots Q^{j_k}$ are admissible
and represent elements of the Dyer-Lashof algebra $\Cal R$, [5], [14].

Denote by $\Cal R_nM$ the submodule of $\Cal R\otimes M$ generated by
the elements  $$Q^{j_1}\dots Q^{j_k}\otimes x$$ where $Q^{j_1}\dots
Q^{j_k},$
the admissible sequences with $dim(x)\le j_k\le dim(x)+n$. Then there is
an isomorphism $\Cal E_nM\cong\Cal R_nM$.

The correspondence $M\longmapsto\Cal R_nM$ determines the monad $\Cal
R_n$ in the category of graded modules.  Algebras over the monad $\Cal
R_n$ we will call $\Cal R_n$-modules.

For an $n$-Lie algebra $L_n$ we will have the module $\Cal R_nL_n$.
Denote also by $T_s\Cal R_nL_n$ the quotient algebra of the free
commutative algebra generated by the module $\Cal R_nL_n$ modulo the
relations $x\cdot x=e_0x$.

From the above considerations it follows that if $X$ is a chain
complex (over $Z/2$) considered as the trivial $E$-coalgebra, then
there are isomorphisms
$$S^{-n}H_*(F(S^nE,E,X))\cong T_s\Cal R_{n-1}L_{n-1}S^{-n}H_*(X).$$
Hence we have

{\bf Theorem 3.} {\sl If $X$ is an $n$-connected topological space
then there is the spectral sequence which converges to $H_*(\Omega^nX)$
and for the first term of this spectral sequence (over $Z/2$) there
is the isomorphism $$E^1\cong T_s\Cal R_{n-1}L_{n-1}S^{-n}H_*(X).$$}

Define the  notion  of  a  $\Cal  P_n$-algebra with $Z/2$-coefficients
generalizing the notion of an $n$-Poisson algebra.

A graded module $M$ (over $Z/2$) will be called a  $\Cal  P_n$-algebra
if

1. There  is  given  a  structure  of a commutative algebra $$x\otimes
y\longmapsto x\cdot y,\quad x,y\in M.$$

2. There is  given  a  structure  of  an  $n$-Lie  algebra  $$x\otimes
y\longmapsto [x,y],\quad x,y\in M,$$ and the $n$-Lie algebra structure
with the commutative algebra structure  form  an  $n$-Poisson  algebra
structure.

3. There is given a structure of $\Cal R_n$-module $$\Cal R_nM\to M,$$
compatible with the $n$-Poisson algebra structure [4].

Denote by $\Cal P_n$ the monad which associates to a graded module $M$
the  free  $\Cal  P_n$-algebra  generated  by  $M$.  Then  there is an
isomorphism $$\Cal P_n(M)\cong T_s\Cal R_nL_nM$$ and the Theorem 3 may
be reformulated

{\bf Theorem  4.}  {\sl  The  first  term  of  the considered spectral
sequence of $H_*(\Omega^nX)$ (over $Z/2$) is isomorphic  to  the  free
$\Cal  P_{n-1}$-algebra generated by $S^{-n}H_*(X)$,  i.e.  $$E^1\cong
\Cal P_{n-1}S^{-n}H_*(X).$$}

Consider now $Z/p$-coefficients, $p>2$. Note that the module $PT_sM$ of
primitive elements of the Hopf algebra $T_sM$ in this case is generated
by the  elements  $x\in  M$ and the elements $x^{p^k}$ for which $k>0$
and $dim(x)$  is even.

However, the  homology  $H_*(FT_sM)$  is  generated  not  only  by the
primitive elements of $T_sM$ but also by the elements $\beta  x^{p^k}$
where  $\beta$  is  the  Bockstein  homomorphism.  Denote  the  module
generated by the primitive elements and the elements  $\beta  x^{p^k}$
by $P_\beta T_sM$.

Then the  homology  of $FT_sM$ will be isomorphic to $T_sS^{-1}P_\beta
T_sM$.

Thus for the homology $H_*(F^nX)$ of the iterated cobar  constructions
$F^nX$    over   $Z/p$-coefficients   there   are   the   isomorphisms
$$H_*(F^nX)\cong T_s(S^{-1}P_\beta T_s)^{n-1}LS^{-1}H_*(X)$$ and hence
the   isomorphisms   $$S^{-n}H_*(F(S^nE,E,X))\cong   T_s(S^{-1}P_\beta
T_s)^{n-1}LS^{-1}H_*(X).$$

Recall that a graded module $L$ (over $Z/p$) is called an $n$-Lie
algebra if there is given an operation $[~,~]\colon L\otimes L\to L$ of
dimension $n$, called an $n$-Lie bracket and satisfying the relations
$$\gather
[x,y]+(-1)^\epsilon [y,x]=0;\\
[x,[y,z]]=[[x,y],z]+(-1)^\epsilon [y,[x,z]];\endgather $$
where $\epsilon =(dim(x)+n)\cdot (dim(y)+n)$.

For a graded module $M$ denote by $\Cal E_n(M)$ the  module  generated
by                            the                            sequences
$\beta^{\epsilon_1}e_{i_1}\dots\beta^{\epsilon_k}e_{i_k}x$,      where
$x\in M$, $\epsilon=0,1$, $0\le i_1\le\dots \le i_k\le n$.

For the elements $e_iy$ it is demanded that $i$ and $dim(y)$ have the
same parity and the dimensions of the elements $e_iy$ are defined equal
to $p\cdot dim(y)+(p-1)i$.

These elements
$\beta^{\epsilon_1}e_{i_1}\dots\beta^{\epsilon_k}e_{i_k}x$,
may be rewritten using the mod-p Dyer-Lashof algebra $\Cal R$
in the form $\beta^{\epsilon_1}Q^{j_1}\dots\beta^{\epsilon_k}
Q^{j_k}x$, where $$\gather 2j_k=i_k+dim(x),\\
2j_{k-1}=i_{k-1}+(p-1)i_k+pdim(x)-\epsilon_k,\\ \dots \\
2j_1=i_1+(p-1)i_2-\epsilon_2+\dots+p^{k-1}dim(x).\endgather $$
It is clear that the sequences $\beta^{\epsilon_1}Q^{j_1}\dots
\beta^{\epsilon_k} Q^{j_k}$ are admissible
and represent elements of the Dyer-Lashof algebra $\Cal R$.

Denote by $\Cal R_nM$ the submodule of $\Cal R\otimes M$ generated by
the elements  $$\beta^{\epsilon_1}Q^{j_1}\dots\beta^{\epsilon_k} Q^{j_k}
\otimes x$$ where $\beta^{\epsilon_1}Q^{j_1}\dots\beta^{\epsilon_k}
Q^{j_k}$
are admissible sequences with $dim(x)\le 2j_k\le dim(x)+n$. Then there
is an isomorphism $\Cal E_nM\cong\Cal R_nM$.

The correspondence  $M\longmapsto\Cal R_nM$ determines the monad $\Cal
R_n$ in the category of graded $Z/p$-modules.  Algebras over the monad
$\Cal R_n$ we will call $\Cal R_n$-modules.

For an $n$-Lie algebra $L_n$ we have the module $\Cal R_nL_n$.
Denote also by $T_s\Cal R_nL_n$ the quotient algebra of the free
commutative algebra generated by the module $\Cal R_nL_n$ modulo the
relations $e_0x=x^p$.

From the above considerations follows that if $X$ is a chain
complex (over $Z/p$) considered as a   trivial $E$-coalgebra, then
there are the isomorphisms
$$S^{-n}H_*(F(S^nE,E,X))\cong T_s\Cal R_{n-1}L_{n-1}S^{-n}X_*.$$
Hence we have

{\bf Theorem 3'.} {\sl For the first term of the considered spectral
sequence of $H_*(\Omega^nX)$ (over $Z/p$), there is an  isomorphism
$$E^1\cong T_s\Cal R_{n-1}L_{n-1}S^{-n}H_*(X).$$}

A graded  module  $M$ (over $Z/p$) will be called a $\Cal P_n$-algebra
if

1.  There is given a structure of a commutative algebra  $$x\otimes
y\longmapsto x\cdot y,\quad x,y\in M.$$

2. There  is  given  a  structure  of  an  $n$-Lie  algebra $$x\otimes
y\longmapsto [x,y],\quad x,y\in M,$$ and the $n$-Lie algebra structure
with  the  commutative  algebra  structure form an $n$-Poisson algebra
structure.

3. There is given a structure of an $\Cal R_n$-module  $$\Cal  R_nM\to
M,$$ compatible with the $n$-Poisson algebra structure [4].

Denote by $\Cal P_n$ the monad which associates to a graded module $M$
the free $\Cal  P_n$-algebra  generated  by  $M$.  Then  there  is  an
isomorphism $$\Cal P_n(M)\cong T_s\Cal R_nL_nM$$ and hence the Theorem
3' may be reformulated

{\bf Theorem 4'.} {\sl The  first  term  of  the  considered  spectral
sequence  of  $H_*(\Omega^nX)$  (over $Z/p$) is isomorphic to the free
$\Cal P_{n-1}$-algebra generated by  $S^{-n}H_*(X)$,  i.e.  $$E^1\cong
\Cal  P_{n-1}S^{-n}H_*(X).$$}

Note that in the case of characteristic zero coefficients the module
$PT_sM$ of primitive elements of $T_sM$ is isomorphic to $M$. Hence
we have

{\bf Theorem 3''.} {\sl For the first term of the considered spectral
sequence of $H_*(\Omega^nX)$ (over a field of characteristic zero)
there is an isomorphism $$E^1\cong T_sL_{n-1}S^{-n}H_*(X).$$}

\vskip .5cm
\centerline{\S 5. The second term of the spectral sequence of the
homology}
\centerline{of iterated loop spaces}
\vskip 6pt
To determine  the second term of the spectral sequence of the homology
of iterated loop spaces, we consider the relation between the Steenrod
algebra and the Dyer-Lashof algebra.

We begin with $Z/2$-coefficients.
Let  $\Cal A$ be the mod-2 Steenrod algebra with $Sq^0\ne 1$ and $\Cal
K$ the Milnor coalgebra (dual to the Steenrod algebra). We will consider
$\Cal K$ as a family $\Cal K(m)$ of polynomial algebras $\Cal K(m)$,
$m\ge 1$ generated by the elements $\xi_i(m)$, $1\le i\le m$ of
dimensions $2^i-1$.

Define products $\Cal K\otimes\Cal K$, $\Cal K\times\Cal K$ to be the
families  $$\gather
(\Cal K\otimes\Cal K)(m)=\Cal K(m)\otimes\Cal K(m),\\
(\Cal K\times\Cal K)(m)=\sum_k\Cal K(k)\otimes\Cal K(m-k).\endgather $$
We denote the multiplication by $\pi\colon\Cal K\otimes\Cal K\to \Cal
K$.

The comultiplication $\nabla\colon\Cal K\to\Cal K\times\Cal K$ is
defined on the generators by the formula
$$\nabla(\xi_i(m))=\sum_{j,k}\xi_{i-j}^{2^j}(k)\otimes\xi_j(m-k).$$

On the other elements the comultiplication is determined by the Hopf
relation. So we can consider $\Cal K$ as a Hopf algebra.

Let now $\Cal R$ be the Dyer-Lashof algebra. Define a homomorphism
$\varphi\colon\Cal K\to\Cal R$ of dimension $-1$ by putting
$$\varphi(x)=\cases Q^{i-1},&x=\xi_1^i(1),\\
0,&~otherwise.\endcases $$
It is easy to see that $\varphi\cup\varphi=0$ and hence $\varphi$ is a
twisting cochain.

Denote by $\overline E_*$ the comonad in the category of graded
modules which associates  to a graded module $M$ the graded
module $H_*(\overline E(M))$. The $E$-coalgebra structure on
the chain complex $C_*(X)$ of a topological space $X$ induces
on its homology $H_*(X)$ the $\overline E_*$-coalgebra structure.

This structure consists of the commutative coalgebra structure
$$\nabla_*\colon H_*(X)\to H_*(X)\otimes H_*(X)$$ and of the coaction
$$\tau\colon H_*(X)\to\Cal K\otimes H_*(X),$$
which are  compatible  in  the  sense  that  $\tau$  is  a  mapping  of
coalgebras.

For an $n$-connected space $X,$ these structures induce on
$\Cal P_{n-1}S^{-n}H_*(X)$ a differential $d_\varphi$ defined on the
generators $s^{-n}x_{i+n}\in H_{i+n}(X)$ by the formula
$$d_\varphi(s^{-n}x_{i+n})=\sum_{x'<x''}[s^{-n}x',s^{-n}x'']+
\varphi\cap s^{-n}x_{i+n}.$$

Denote the corresponding differential $\Cal P_{n-1}$-algebra by
$\Cal P_{n-1\varphi}S^{-n}H_*(X)$. Then we will have

{\bf Theorem  5.} {\sl For the second term of the spectral sequence of
$H_*(\Omega^nX)$ (over  $Z/2$)  there  is  an  isomorphism  $$E^2\cong
H_*(\Cal P_{n-1\varphi}S^{-n}H_*(X)).$$}

Consider now the relation between the mod-p Steenrod algebra and the mod-p 
Dyer-Lashof algebra. 
Let  $\Cal A$ be the mod-p Steenrod algebra with $Sq^0\ne 1$ and $\Cal
K$ the Milnor coalgebra (dual to the Steenrod algebra). We will consider
$\Cal K$ as a family $\Cal K(m)$ of polynomial algebras $\Cal K(m)$,
$m\ge 1$ generated by the elements $\xi_i(m)$, $\tau_i(m)$,
$1\le i\le m$ of dimensions $2(p^i-1)$ and $2p^i-1$ correspondingly.

There are a multiplication $$\pi\colon\Cal K\otimes\Cal K\to\Cal
K$$ and a comultiplication $$\nabla\colon\Cal K\to\Cal K\times\Cal K$$
defined     on      the      generators      by      the      formulas
$$\nabla(\xi_i(m))=\sum_{j,k}\xi_{i-j}^{p^j}(k)\otimes\xi_j(m-k).$$
$$\nabla(\tau_i(m))=\tau_i(m)\otimes                              1+
\sum_{j,k}\xi_{i-j}^{p^j}(k)\otimes\tau_j(m-k).$$    On    the   other
elements the comultiplication is determined by the Hopf  relation.  So
we can consider $\Cal K$ as a Hopf algebra.

Let now $\Cal R$ be the mod-p Dyer-Lashof algebra. Define the
homomorphism
$\varphi\colon\Cal K\to\Cal R$ of dimension $-1$ by putting
$$\varphi(x)=\cases \beta P^i,&x=\xi_1^i(1),\\
P^i,&x=\xi_1^{i-1}(1)\tau_1(1),\\
0,&~otherwise.\endcases $$

It is easy to see that $\varphi\cup\varphi=0$ and hence $\varphi$ is
the twisting cochain.

Denote by $\overline E_*$ the comonad in the category of graded
$Z/p$-modules which associates  to a graded module $M$ the graded
module $H_*(\overline E(M))$. The $E$-coalgebra structure on
the chain complex $C_*(X)$ of a topological space $X$ induces
on its homology $H_*(X)$ an  $\overline E_*$-coalgebra structure.

This structure consists of the commutative coalgebra structure
$$\nabla_*\colon H_*(X)\to H_*(X)\otimes H_*(X)$$ and the
coaction of the Milnor coalgebra
$$H_*(X)\to\Cal K\otimes H_*(X).$$

For an $n$-connected space $X,$ these structures induce on
$\Cal P_{n-1}S^{-n}H_*(X)$ a differential $d_\varphi$ defined on the
generators $s^{-n}x_{i+n}\in H_{i+n}(X)$ by the formula
$$d_\varphi(s^{-n}x_{i+n})=\sum_{x'<x''}(-1)^\epsilon
[s^{-n}x',s^{-n}x'']+
\varphi\cap s^{-n}x_{i+n},$$ where $\sum x'\otimes
x''=\nabla_*(x_{i+n})$,
$\epsilon=n\cdot dim(x')$.

Denote the corresponding differential $\Cal P_{n-1}$-algebra by
$\Cal P_{n-1\varphi}S^{-n}H_*(X)$.
Then we have

{\bf Theorem 5'.} {\sl For the second term of the considered
spectral sequence of $H_*(\Omega^nX)$ (over $Z/p$) there is an
isomorphism $$E^2\cong H_*(\Cal P_{n-1\varphi}S^{l-n}H_*(X)).$$}

Note that in the case of characteristic zero coefficients, the
$\overline E_*$-coalgebra structure on the homology $H_*(X)$
of a topological space $X$ consists only of the commutative
coalgebra structure
$$\nabla_*\colon H_*(X)\to H_*(X)\otimes H_*(X).$$
Hence we have

{\bf Theorem 5''.} {\sl For the second term of the considered spectral
sequence of $H_*(\Omega^nX)$ (over a field of characteristic zero),
there is an  isomorphism
$$E^2\cong H_*(P_{n-1\varphi}S^{-n}H_*(X)),$$
where the differential $d_\varphi$ is defined on the generators
$s^{-n}x_{i+n}\in H_{i+n}(X)$ by the formula
$$d_\varphi(s^{-n}x_{i+n})=\sum_{x'<x''}(-1)^\epsilon
[s^{-n}x',s^{-n}x''],$$
where $\sum x'\otimes x''=\nabla_*(x_{i+n})$, $\epsilon=n\cdot
dim(x')$.}

\vskip .5cm
\centerline{\S 6. The homology of iterated loop spaces of   the
real projective spaces}
\vskip 6pt
Some of the most important spaces in Algebraic Topology, besides
spheres, are the real projective spaces $RP^n$ and $RP^\infty$.
In some sense these spaces are opposite to spheres. The homology
$H_*(RP^\infty)$ (over $Z/2$) is the free $\overline E_*$-coalgebra
with one $1$-dimensional generator, whereas for spheres that homology
is a trivial $\overline E_*$-coalgebra.

On   the other hand, the Adams spectral sequences of stable homotopy
groups of  $RP^{\infty}$  and $RP^\infty/RP^n$ is very similar to the
corresponding spectral sequence of spheres;  it is a problem  to  find
the relations between them.

The homology  of  iterated  loop spaces gives an approximation of the
homotopy groups.  Here we consider the problem of calculation  of  the
homology   $$H_*(\Omega^m(RP^\infty/RP^n)),\quad   m\le  n,$$  of  the
iterated  loop  spaces  of  the  space  $RP^\infty/RP^n$  with   $Z/2$
coefficients.

Denote the $i$-dimensional generator of $H_*(RP^\infty)$ as $e_i$.
A coalgebra structure on these generators is determined by the
formula $$\nabla(e_i)=\sum_je_j\otimes e_{i-j}.$$

An action of the Steenrod algebra is determined by the formula
$$Sq^j(e_i)=\binom{i-j}je_{i-j}.$$

Usually, on the homology $H_*(X)$ of a topological space $X$ (over $Z/2$)
there is not only $\overline E_*$-coalgebra structure. Besides that there 
are functional homology operations $$\tau_*^n\colon H_*(X)\to\overline
E_*^n(H_*(X)),\quad [15], [16].$$
These operations determine the higher differentials of the spectral
sequence.

Fortunately, $H_*(RP^\infty)$ is a free $\overline  E_*$-coalgebra  and
hence  the higher functional  homology  operations  are  trivial.  
The  homology $H_*(RP^n)$ is the $\overline  E_*$-subcoalgebra  of  
$H_*(RP^\infty)$ and  hence the $\overline E_*$-coalgebra structure on 
$H_*(RP^\infty)$ induces an $\overline E_*$-coalgebra structure   on
$$H_*(RP^\infty)/H_*(RP^n)\cong  H_*(RP^\infty/RP^n).$$  Therefore  on
the homology $H_*(RP^\infty/RP^n)$  there  are  no  higher  functional
homology operations.
So the spectral sequence of the homology of the iterated loop spaces of
$RP^\infty/RP^n$ has no higher differentials and we have

{\bf Theorem 6.} {\sl The homology $H_*(\Omega^m(RP^\infty/RP^n))$,
$n\ge m$ (over $Z/2$) is isomorphic to the homology of the differential
$\Cal P_{m-1}$-algebra $\Cal P_{m-1\varphi}S^{-m}H_*(RP^\infty/RP^n)$,
where the differential $d_\varphi$ on the generators
$s^{-m}e_i$ is defined by the formula
$$d_\varphi(s^{-m}e_i)=\sum_{\scriptstyle j\atop \scriptstyle 2j<i}
[s^{-m}e_j,s^{-m}e_{i-j}]+\sum_{\scriptstyle j\atop \scriptstyle 2j\le
i}\binom{i-j}je_{2j-i+m-1}(s^{-m}e_{i-j}).$$}

Note that in the case $m=1$ we have the isomorphisms
$$\Cal P_0S^{-1}H_*(RP^\infty/RP^n)\cong T_sLS^{-1}H_*(RP^\infty/
RP^n)\cong TS^{-1}H_*(RP^\infty/RP^n),$$
and the differential in the tensor algebra
$TS^{-1}H_*(RP^\infty/RP^n)$ on the generators $u_i=s^{-1}e_{i+1}$ has
a very simple form:
$$d_\varphi(u_i)=\sum_ju_j\otimes u_{i-j-1}.$$ From here it follows that
the homology $H_*(\Omega(RP^\infty/RP^n))$ (over $Z/2$) is isomorphic
to the algebra generated by the elements $u_i$, $n\le i\le 2n$,
of dimensions $i$ and relations
$$\gather u_n\cdot u_n=0;\\
u_n\cdot u_{n+1}+u_{n+1}\cdot u_n=0;\\
\dots \\
u_n\cdot u_{2n}+u_{n+1}\cdot u_{2n-1}+\dots+u_{2n}\cdot u_n=0.
\endgather $$

So, as a graded module, the homology $H_*(\Omega(RP^\infty/RP^n))$
is generated by the noncommutative products $u_{n_1}\cdot\dots
\cdot u_{n_k}$ with $n\le n_1\le 2n$, $n<n_2,\dots,n_k\le 2n$.

Consider the homology $H_*(\Omega^2(RP^\infty/RP^2))$. It
is isomorphic to the homology of the differential $\Cal P_1$-algebra
$\Cal P_{1\varphi}\{u_i|~i\ge 1\}$, where $u_i=s^{-2}e_{i+2}$
and the differential is determined by the formulas
$$\align
d_\varphi(u_{2i+1})&=\sum_{j<i}[u_j,u_{2i-j-1}]+\binom{i+2}{i+1}
e_0(u_i);\\
d_\varphi(u_{2i+2})&=\sum_{j<i}[u_j,u_{2i-j}]+e_1(u_i).
\endalign $$

In small dimensions we have

$d(u_1)=0$;

$d(u_2)=0$;

$d(u_3)=u_1u_1$;

$d(u_4)=e_1(u_1)$;

$d(u_5)=[u_1,u_2]$;

$d(u_6)=[u_1,u_3]+e_1(u_2)$;

$d(u_7)=[u_1,u_4]+[u_2,u_3]+u_3u_3$;

$d(u_8)=[u_1,u_5]+[u_2,u_4]+e_1(u_3)$.

\vskip 6pt
From these formulas it follows that the homology $H_i=H_i(\Omega^2
(RP^\infty/RP^2))$ in small dimensions $i$ has the following
generators

$H_1:\quad u_1$;

$H_2:\quad u_2$;

$H_3:\quad u_1u_2$;

$H_4:\quad u_2^2$;

$H_5:\quad u_1u_2^2,~ e_1(u_2)$;

$H_6:\quad u_1e_1(u_2),~ u_2^3,~ [u_2,u_3],~ u_3^2$;

$H_7:\quad u_1u_2^3,~ u_1u_3^2,~ u_1[u_2,u_3],~ u_2e_1(u_2),~
e_1(u_3)$.

\vskip 6pt

Consider the homology $H_*(\Omega^3(RP^\infty/RP^3))$. It
is isomorphic to the homology of the differential $\Cal P_2$-algebra
$\Cal P_{2\varphi}\{u_i|~i\ge 1\}$, where $u_i=s^{-3}e_{i+3}$
and the differential is determined by the formulas
$$\align
d_\varphi(u_{2i+1})&=\sum_{j<i-1}[u_j,u_{2i-j-2}]+\binom{i+3}{i+1}
e_0(u_i)+e_2(u_{i-1});\\
d_\varphi(u_{2i+2})&=\sum_{j\le i-1}[u_j,u_{2i-j-1}]+
\binom{i+3}{i+2}e_1(u_i).\endalign $$

In small dimensions we  have

$d(u_1)=0$;

$d(u_2)=0$;

$d(u_3)=0$;

$d(u_4)=0$;

$d(u_5)=e_2(u_1)$;

$d(u_6)=[u_1,u_2]+e_1(u_2)$;

$d(u_7)=[u_1,u_3]+e_2(u_2)+u_3u_3$;

$d(u_8)=[u_1,u_4]+[u_2,u_3]$.

\vskip 6pt

From these formulas it follows that the homology $H_i=H_i(\Omega^3
(RP^\infty/RP^3))$ in small dimensions $i$ has the following
generators

$H_1:\quad u_1$;

$H_2:\quad u_1^2,~ u_2$;

$H_3:\quad u_1^3,~ u_1u_2,~ e_1(u_1),~ u_3$;

$H_4:\quad u_1^4,~ u_1^2u_2,~ u_1e_1(u_1),~ u_1u_3,~ u_2^2,~ u_4$;

$H_5:\quad u_1^5,~u_1^3u_2,~u_1^2u_3,~u_1^2e_1(u_1),~u_1u_2^2,~
u_1u_4,~u_2u_3,~u_2e_1(u_1),~e_1(u_2)$;

$H_6:\quad u_1^6,~u_1^4u_2,~u_1^3u_3,~u_1^3e_1(u_1),~u_1^2u_2^2,~
u_1^2u_4,~u_1u_2u_3,~u_1u_2e_1(u_1),~u_1e_1(u_2),$
$u_2u_4,~u_2^3,~u_3e_1(u_1),~u_3^2,~e_2(u_2),~e_1(u_1)^2$;

$H_7:\quad u_1^7,~u_1^5u_2,~u_1^4u_3,~u_1^4e_1(u_1),~u_1^3u_2^2,~
u_1^3u_4,~u_1^2u_2u_3,~u_1^2u_2e_1(u_1),~u_1^2e_1(u_2),$
$u_1u_2u_4,~u_1u_2^3,~u_1u_3e_1(u_1),~u_1u_3^2,~u_1e_2(u_2),
~u_1e_1(u_1)^2,~u_2^2u_3,~u_2^2e_1(u_1)$,$u_2e_1(u_2),$
$u_3u_4,~[u_2,u_3],~e_1(u_3),~e_1e_1(u_1)$.

\vskip 6pt
Consider also the homology $H_*(\Omega^4(RP^\infty/RP^4))$. It
is isomorphic to the homology of the differential $\Cal P_3$-algebra
$\Cal P_{3\varphi}\{u_i|~i\ge 1\}$, where $u_i=s^{-4}e_{i+4}$
and the differential is determined by the formulas
$$\align
d_\varphi(u_{2i+1})&=\sum_{j<i-1}[u_j,u_{2i-j-3}]+\binom{i+4}{i+1}
e_0(u_i)+\binom{i+3}{i+2}e_2(u_{i-1});\\
d_\varphi(u_{2i+2})&=\sum_{j<i-1}[u_j,u_{2i-j-2}]+
\binom{i+4}{i+2}e_1(u_i)+e_3(u_{i-1}).\endalign $$

In small dimensions we have

$d(u_1)=0$;

$d(u_2)=0$;

$d(u_3)=0$;

$d(u_4)=0$;

$d(u_5)=e_2(u_1)$;

$d(u_6)=e_1(u_2)+e_3(u_1)$;

$d(u_7)=[u_1,u_2]+u_3u_3$;

$d(u_8)=[u_1,u_3]+e_1(u_3)+e_3(u_2)$.

\vskip 6pt

From these formulas it follows that the homology $H_i=H_i(\Omega^4
(RP^\infty/RP^4))$ in small dimensions $i$ has the following
generators

$H_1:\quad u_1$;

$H_2:\quad u_1^2,~ u_2$;

$H_3:\quad u_1^3,~ u_1u_2,~ e_1(u_1),~ u_3$;

$H_4:\quad u_1^4,~ u_1^2u_2,~ u_1e_1(u_1),~ u_1u_3,~ u_2^2,~ u_4$;

$H_5:\quad u_1^5,~u_1^3u_2,~u_1^2u_3,~u_1^2e_1(u_1),~u_1u_2^2,~
u_1u_4,~u_2u_3,~u_2e_1(u_1),~e_1(u_2)$;

$H_6:\quad u_1^6,~u_1^4u_2,~u_1^3u_3,~u_1^3e_1(u_1),~u_1^2u_2^2,~
u_1^2u_4,~u_1u_2u_3,~u_1u_2e_1(u_1),~u_1e_1(u_2)$,
$u_2^3,~u_2u_4,~u_3e_1(u_1),~u_3^2,~e_2(u_2),~e_1(u_1)^2$;

$H_7:\quad u_1^7,~u_1^5u_2,~u_1^4u_3,~u_1^4e_1(u_1),~u_1^3u_2^2,~
u_1^3u_4,~u_1^2u_2u_3,~u_1^2u_2e_1(u_1),~u_1^2e_1(u_2)$,
$u_1u_2^3,~u_1u_2u_4,~u_1u_3e_1(u_1),~u_1u_3^2,~u_1e_2(u_2),~
u_1e_1(u_1)^2,~u_2^2u_3,~u_2^2e_1(u_1),~u_2e_1(u_2)$,
$u_3u_4,~e_1(u_3),~e_3(u_2),~e_1e_1(u_1)$.

\vskip .5cm
\centerline{\S 7. The homology of iterated loop spaces of   the
complex projective spaces}
\vskip 6pt
Here we consider the problem of calculation the homology
$$H_*(\Omega^m(CP^\infty/CP^n)),\quad m\le 2n+1,$$ of iterated loop 
spaces of the space $CP^\infty/CP^n$. We begin with $Z/2$ coefficients.

Denote the $2i$-dimensional generator of $H_*(CP^\infty)$ as $c_i$.
A coalgebra structure on this generators is determined by the formula
$$\nabla(c_i)=\sum_jc_j\otimes c_{i-j}.$$

An action of the Steenrod algebra is determined by the formula
$$Sq^{2j}(c_i)=\binom{i-j}jc_{i-j}.$$

As above, since $H_*(CP^\infty)$ is the free $\overline E_*$-coalgebra 
there are no higher $\overline E_*$-operations.  The homology $H_*(CP^n)$ 
is the $\overline   E_*$-subcoalgebra   of  $H_*(CP^\infty)$  and  hence  
the $\overline E_*$-coalgebra structure on  $H_*(CP^\infty)$  induces  the
$\overline       E_*$-coalgebra       structure      on      $$H_*(CP^
\infty)/H_*(CP^n)\cong H_*(CP^\infty/CP^n)$$ Therefore on the homology
$H_*(CP^\infty/CP^n)$  there are no higher $\overline E_*$-operations.
So the spectral sequence of the homology of the iterated  loop  spaces
of $CP^\infty/CP^n$ has no higher differentials and we have

{\bf Theorem 7.} {\sl The homology $H_*(\Omega^m(CP^\infty/CP^n))$,
$m\le 2n+1$ (over $Z/2$) is isomorphic to the homology of the
$\Cal P_{m-1}$-algebra $\Cal P_{m-1\varphi}S^{-m}
H_*(CP^\infty/CP^n)$, where the differential $d_\varphi$ on the
generators $s^{-m}c_i$ is defined by the formula
$$d_\varphi(s^{-m}c_i)=\sum_{\scriptstyle j\atop \scriptstyle 2j<i}
[s^{-m}c_j,s^{-m}c_{i-j}]+\sum_{\scriptstyle j\atop\scriptstyle 2j\le i}
\binom{i-j}je_{2(2j-i)+m-1}(s^{-m}c_{i-j}).$$}

Note that in the case $m=1$ we have the isomorphisms
$$\Cal P_0S^{-1}H_*(CP^\infty/CP^n)\cong T_sLS^{-1}H_*(CP^\infty/
CP^n)\cong TS^{-1}H_*(CP^\infty/CP^n),$$
and the differential in the tensor algebra
$TS^{-1}H_*(CP^\infty/CP^n)$ on the generators $v_i=s^{-1}c_{i+1}$
of dimensions $2i+1$ has a very simple form:
$$d_\varphi(v_i)=\sum_jv_j\otimes v_{i-j-1}.$$

From here it follows that the homology $H_*(\Omega(CP^\infty/CP^n))$
(over $Z/2$) is isomorphic to the algebra generated by the elements
$v_i$, $n\le i\le 2n$, and relations
$$\gather v_n\cdot v_n=0;\\
v_n\cdot v_{n+1}+v_{n+1}\cdot v_n=0;\\
\dots \\
v_n\cdot v_{2n}+v_{n+1}\cdot v_{2n-1}+\dots+v_{2n}\cdot v_n=0.\endgather $$

So as a graded module the homology $H_*(\Omega(CP^\infty/CP^n))$
is generated by the noncommutative products $v_{n_1}\cdot\dots
\cdot v_{n_k}$ with $n\le n_1\le 2n$, $n<n_2,\dots,n_k\le 2n$.

Consider the   homology   $H_*(\Omega^2(CP^\infty/CP^n))$.    It    is
isomorphic  to  the  homology  of  the differential $\Cal P_1$-algebra
$\Cal   P_{1\varphi}\{v_i|~i\ge   n\}$,   where   $v_i=s^{-2}c_{i+1}$,
$dim(v_i)=2i$  and  the  differential  is determined by the formulas
$$\align            d_\varphi(v_{2i})&=\sum_{j<i}[v_j,v_{2i-j-1}];\\
d_\varphi(v_{2i+1})&=\sum_{j<i}[v_j,v_{2i-j}]+e_1(v_i).\endalign $$

From these formulas it follows that the homology $H_*(\Omega^2
(CP^\infty/CP^n))$ is the $\Cal P_1$-algebra generated by the elements
$v_n,\dots ,v_{2n}$ and relations
$$\gather e_1(v_n)=0;\\ 
[v_n,v_{n+1}]=0;\\[v_n,v_{n+2}]+e_1(v_{n+1})=0;\\
\dots \\ [v_n,v_{2n}]+[v_{n+1},v_{2n-1}]+\dots =0.\endgather $$

Consider the   homology   $H_*(\Omega^4(CP^\infty/CP^2))$.    It    is
isomorphic  to  the  homology  of  the differential $\Cal P_3$-algebra
$\Cal  P_{3\varphi}\{v_i|~i\ge  1\}$,   where   $v_i=s^{-4}c_{i+2}$,
$dim(v_i)=2i$  and  the  differential  is determined by the formulas
$$\align
d_\varphi(v_{2i+1})&=\sum_{j<i}[v_j,v_{2i-j-1}]+
\binom{i+2}{i+1}e_1(v_i);\\
d_\varphi(v_{2i+2})&=\sum_{j<i}[v_j,v_{2i-j}]+e_3(v_i).\endalign $$

In small dimensions we have

$d(v_1)=0$;

$d(v_2)=0$;

$d(v_3)=e_1(v_1)$;

$d(v_4)=e_3(v_1)$;

$d(v_5)=[v_1,v_2]$;

$d(v_6)=[v_1,v_3]+e_3(v_2)$;

$d(v_7)=[v_1,v_4]+[v_2,v_3]+e_1(v_3)$;

\vskip 6pt

From these formulas it follows that the homology $H_i=H_i(\Omega^4
(CP^\infty/CP^2))$ in small dimensions $i$ has the following
generators

$H_1=0$;

$H_2:\quad v_1$;

$H_3=0$;

$H_4:\quad v_1^2,~ v_2$;

$H_5=0$;

$H_6:\quad v_1^3,~ v_1v_2,~ e_2(v_1)$;

$H_7=0$;

$H_8:\quad v_1^4,~ v_1^2v_2,~ v_1e_2(v_1),~ v_2^2$;

$H_9:\quad e_1(v_2)$;

$H_{10}:\quad v_1^5,~v_1^3v_2,~v_1^2e_2(v_1),~v_1v_2^2,~
v_2e_2(v_1),~e_2(v_2)$;

$H_{11}:\quad v_1e_1(v_2),~e_3(v_2)$;

$H_{12}:\quad v_1^6,~v_1^4v_2,~v_1^3e_2(v_1),~v_1^2v_2^2,~
v_1v_2e_2(v_1),~v_1e_2(v_2)$;

$H_{13}:\quad v_1^2e_1(v_2),~v_1e_3(v_2),~e_1(v_3),~[v_2,v_3]$.

\vskip 6pt
Consider now $Z/p$ coefficients, $p>2$.
The coalgebra structure on the generators is determined by the
formula
$$\nabla(c_i)=\sum_jc_j\otimes c_{i-j}.$$

The action of the mod-p Steenrod algebra is determined by the
formula
$$P^j(c_i)=\binom{i-(p-1)j}jc_{i-(p-1)j}.$$

Since $H_*(CP^\infty)$  is  the free $\overline E_*$-coalgebra,  there
are no higher $\overline E_*$-operations.  The homology $H_*(CP^n)$ is
the  $\overline  E_*$-subcoalgebra  of  $H_*(CP^\infty)$ and hence the
$\overline E_*$-coalgebra structure on  $H_*(CP^\infty)$  induces  the
$\overline       E_*$-coalgebra       structure      on      $$H_*(CP^
\infty)/H_*(CP^n)\cong   H_*(CP^\infty/CP^n).$$   Therefore   on   the
homology   $H_*(CP^\infty/CP^n)$   there   are  no  higher  $\overline
E_*$-operations.  So the spectral sequence  of  the  homology  of  the
iterated  loop  spaces of $CP^\infty/CP^n$ has no higher differentials
and we have

{\bf Theorem 8.} {\sl The homology $H_*(\Omega^m(CP^\infty/CP^n))$,
$m\le 2n+1$ (over $Z/p$) is isomorphic to the homology of the
$\Cal P_{m-1}$-algebra $\Cal P_{m-1\varphi}S^{-m}H_*(CP^\infty/CP^n)$,
where the differential $d_\varphi$ on the generators
$s^{-m}e_i$ is defined by the formula
$$\align d_\varphi(s^{-m}c_i)=&
\sum_{\scriptstyle j\atop \scriptstyle 2j\le i}
[s^{-m}c_j,s^{-m}c_{i-j}]+\\+&\sum_{\scriptstyle j\atop \scriptstyle
pj<i}
\binom{i-(p-1)j}j\beta e_{2(pj-i)+m}(s^{-m}c_{i-(p-1)j}).\endalign$$}

Note that in the case $m=1$ we have the isomorphisms
$$\Cal P_0S^{-1}H_*(CP^\infty/CP^n)\cong T_sLS^{-1}H_*(CP^\infty/
CP^n)\cong TS^{-1}H_*(CP^\infty/CP^n),$$
and the differential in the tensor algebra
$TS^{-1}H_*(CP^\infty/CP^n)$ has very simple form
$$d_\varphi(s^{-1}c_i)=\sum_js^{-1}c_j\otimes s^{-1}c_{i-j}.$$

From here it follows that the homology $H_*(\Omega(CP^\infty/CP^n))$
(over $Z/p$) is isomorphic to the algebra generated by the elements
$v_i=s^{-1}c_{i+1}$, $n\le i\le 2n$, of dimensions $2i+1$ and
relations
$$\gather v_n\cdot v_n=0;\\
v_n\cdot v_{n+1}+v_{n+1}\cdot v_n=0;\\
\dots \\
v_n\cdot v_{2n}+v_{n+1}\cdot v_{2n-1}+\dots+v_{2n}\cdot v_n=0.\endgather $$

So as a graded module the homology $H_*(\Omega(CP^\infty/CP^n))$
is generated by the noncommutative products $v_{n_1}\cdot\dots
\cdot v_{n_k}$ with $n\le n_1\le 2n$, $n<n_2,\dots,n_k\le 2n$.

Consider the    homology    $H_*(\Omega^2(CP^\infty/CP^n))$.   It   is
isomorphic to the homology  of  the  differential  $\Cal  P_1$-algebra
$\Cal P_{1\varphi}\{v_i|~i\ge 1\}$,  where $v_i=s^{-2}c_{i+1}$,  $i\ge
n$  and  the  differential  is   determined   by   the   formulas   $$
d_\varphi(v_i)=\sum_{2j<i}[v_j,v_{i-j-1}].$$

From these formulas it follows that the homology $H_i=H_i(\Omega^2
(CP^\infty/CP^n))$ is the $\Cal P_1$-algebra generated by the elements
$v_n,\dots,v_{2n}$ and relations
$$\gather [v_n,v_n]=0;\\ [v_n,v_{n+1}]=0;\\
[v_n,v_{n+2}]+[v_{n+1},v_{n+1}]=0;\\
\dots \\ [v_n,v_{2n}]+[v_{n+1},v_{2n-1}]+\dots =0.\endgather $$

Consider the   homology   $H_*(\Omega^4(CP^\infty/CP^2))$.    It    is
isomorphic  to  the  homology  of  the differential $\Cal P_3$-algebra
$\Cal  P_{3\varphi}\{v_i|~i\ge  1\}$,   where   $v_i=s^{-4}c_{i+2}$,
$dim(v_i)=2i$  and  the  differential  is determined by the formulas
$$\align        d_\varphi(v_{pi-1})&=\sum_{2j<pi-2}[v_j,v_{pi-j-3}]+
\binom{i+1}i\beta                                      e_2(v_{i-1});\\
d_\varphi(v_k)&=\sum_{2j<k-1}[v_j,v_{k-j-2}],~k\ne pi-1.\endalign $$

In small dimensions, we have

$d(v_1)=0$;

$d(v_2)=0$;

$d(v_3)=0$;

$d(v_4)=[v_1,v_1]$;

$d(v_5)=[v_1,v_2]$;

$d(v_6)=[v_1,v_3]+[v_2,v_2]$;

$\dots $

$d(v_{2p-1})=[v_1,v_{2p-4}]+\dots+[v_{p-2},v_{p-1}]+3\beta e_2(v_1)$;

$\dots$

\vskip 6pt
From these formulas it follows that the homology $H_i=H_i(\Omega^4
(CP^\infty/CP^2))$ with, for example, $Z/3$-coefficients, in small
dimensions $i$ has the following generators

$H_1=0$;

$H_2:\quad v_1$;

$H_3=0$;

$H_4:\quad v_1^2,~ v_2$;

$H_5=0$;

$H_6:\quad v_1^3,~ v_1v_2,~ v_3$;

$H_7=0$;

$H_8:\quad v_1^4,~ v_1^2v_2,~ v_1v_3,~ v_2^2$;

$H_9:\quad \beta e_2(v_1)$;

$H_{10}:\quad v_1^5,~v_1^3v_2,~v_1^2v_3,~v_1v_2^2,~v_2v_3,~e_2(v_1)$;

$H_{11}:\quad v_1\beta e_2(v_1),~[v_2,v_2]$.

\vskip 1cm

\centerline{APPENDIX}
\vskip 6pt
\centerline{\sl F. Sergeraert}
\vskip .5cm

The general ideas of the paper [11] led Julio Rubio and myself
to a simple method solving the computability problem for the homology
groups of iterated loop spaces, when the initial space is sufficiently
reduced. The main ingredient is {\it functional Programming}, allowing
us to  constructively  apply  three  particular  cases  of  the  basic
homological perturbation lemma in situations involving highly infinite
simplicial sets; see [23].

This is not only a theoretical result. The algorithm, the existence of
which has been so proved, has been concretely written in Common Lisp
with significant results. The first version of our program, named
{\it EAT} (Effective Algebraic Topology), is Web-reachable
[24] with a rich documentation (250pp.). The most recent version
of this program named {\it Kenzo}\footnote""{Kenzo is the name of my
{\it cat}, and CAT = {\it Constructive Algebraic Topology}.} had just
been finished when this appendix was written (January 1999). See
[25] for a few explanations and a small demonstration file. A
public version of this program, with a reasonably complete user guide,
will be soon distributed at the same Web-address.

Let us give a typical example of the use of Kenzo, related to the main
subject of this paper. What about the first $\Bbb Z$-homology groups of
the iterated loop spaces $\Omega^p_n =\Omega^p (RP^\infty/RP^n)$,
if $p\leq n$? The situation is difficult: the projective space is not a
suspension and the classical results about the homology groups of spaces
of the form $\Omega^nS^nX$ cannot be applied; see [18] for a
survey of this subject. The only general solution known at this time to
determine the groups $H_n(\Omega_n^p;\Bbb Z_2)$, due to Vladimir
Smirnov, is  the  subject  of the main part of this paper.  We explain
here how the results of Smirnov can be verified with the Kenzo program
for the first homology groups.

Let us consider for example the computation of $H_5(\Omega^2 (RP^\infty
/RP^2))$. In an interactive Lisp session where the Kenzo
program has been loaded, you can execute:

{\tt USER(1): (setf trunc3-proj-space (r-proj-space 3)) ==>

[K1 Simplicial-Set]}

This Lisp dialog must be understood as follows. The string
{\tt USER(1):} is the {\it Lisp prompt}; Lisp is waiting
for the next expression you want to evaluate. In this case
the instruction is :

{\tt (setf trunc3-proj-space (r-proj-space 3))}

\noindent   which  means   you   want  to   assign   to  the   symbol
{\tt trunc3-proj-space}   the   result    of   the   evaluation   of
{\tt (r-proj-space 3)}. The last evaluation {\it constructs} in the
Lisp  environment  our  version  of  the  truncated  projective  space
$RP^\infty/RP^2$, and  the resulting  object is assigned  to
the symbol;  furthermore a simple external form of the assigned object
is  displayed;  this  form  can be read:  this is the Kenzo object {\tt 
\#1 K1},  which is a simplicial set. Of course the internal form is
much  more  complicated:  the  internal object {\it is} the looked-for
space, or more precisely {\it codes} this space.

But  we  want to  consider  the  second  loop-space of  the  truncated
projective space. The Kenzo program constructs it in this way:

{\tt USER(2): (setf omega2-trunc3-proj-space

$\phantom{...........................}$(loop-space trunc3-proj-space 2)) ==>

[K22 Simplicial-Group]}

The  result, the  Kenzo  object \# 22, is  a simplicial  {\it group},
namely the  Kan model  of the second  loop space [21],  a highly
infinite object.

The homology groups of an object such as a chain complex, a simplicial
set,  are computed by the Kenzo function {\tt homology}. For example
if  you  want  to verify that the truncated projective space has the
right homology in dimension $5$:

{\tt USER(3): (homology trunc3-proj-space 5) ==>

Computing boundary-matrix in dimension 5.

Rank of the source-module : 1.

Computing boundary-matrix in dimension 6.

Rank of the source-module : 1.

Homology in dimension 5 :

Component Z/2Z

---done--- }

The  homology  group  is  $Z_2$.  The  homology  group  $H_5(\Omega^2
(RP^\infty/RP^2))$ is obtained in the same way:

{\tt USER(4): (homology omega2-trunc3-proj-space 5) ==>

Computing boundary-matrix in dimension 5.

Rank of the source-module : 29.

Computing boundary-matrix in dimension 6.

Rank of the source-module : 70.

Homology in dimension 5 :

Component Z/2Z

Component Z/2Z

---done---}

So that $H_5(\Omega^2 (RP^\infty/RP^2))  = Z_2^2$.
As you see, the program determines a  {\it Hirsch $Z$-complex}
for our loop space,
with  29 generators  in dimension  5  and 70  generators in  dimension
6. Note   the   program   gives   you   the   $Z$-homology,   not   the
$Z_2$-homology, as   in Smirnov's  text. The complexity is very high:
the highest group computed by Kenzo for this second loop space is $H_7
= Z_2^2 \oplus Z_8$; one  day of  CPU time  has been  needed  on a
powerful Linux PC to obtain it.

The topologists frequently ask  for some homology groups determined by
the Kenzo program that are not reachable by human computers. The first
example of  this sort was given  by the previous version  of the Kenzo
program, the EAT program [24]. We repeat it here, because it is
very simple and ten years after its first computation, we have not yet
found a topologist knowing how  to compute this group by hand.  Please
try to do it!

The construction  is the following.  The homotopy  group $\pi_2 \Omega
S^3$ is $Z$, so that attaching a 3-cell to  $\Omega S^3$ by a  map of
degree $2$  makes   sense.   Let   $D   \Omega  S^3$ the space so obtained.
Question: what about the homology groups of $\Omega D\Omega S^3$?
Let us  show how the first homology groups  of this strange loop
space  are computed  by  the  Kenzo program.  Firstly  the loop  space
$\Omega S^3$ is constructed as before:

{\tt USER(5): (setf s3 (sphere 3))  ==>

[K262 Simplicial-Set]

USER(6): (setf os3 (loop-space s3)) ==>

[K267 Simplicial-Group]}

How to attach a $3$-cell?  This $3$-cell will be a $3$-simplex and, in
order  to  attach  it,  we must describe what $2$-simplices of $\Omega
S^3$ are its faces. The list of faces is defined and used as follows:

{\tt USER(7): (setf faces (list (loop3 0 's3 1)

$\phantom{.........................................}$(absm 3 +null-loop+)

$\phantom{.........................................}$(loop3 0 's3 1)

$\phantom{.........................................}$(absm 3 +null-loop+)))
==>

(<<Loop[S3]>> <AbSm 1-XL:Loop>>> <<Loop[S3]>> <AbSm 1-0 <<Loop>>>)

USER(8): (setf dos3 (disk-pasting os3 3 'new faces))

[K380 Simplicial-Set]}

The  {\tt disk-pasting}  function  atttaches  a $3$-cell  of  ``name''
{\tt new} to  {\tt os3}, using  the face list  {\tt faces}; this
list is:

\item 0) The ``fundamental'' simplex of $\Omega S^3$;

\item 1) The second face is collapsed on the base point;

\item 2) The same as 0);

\item 3) The same as 2);

\noindent   so  that, taking  account  of  the  usual sign  rules, the
attaching map has  degree $2$.  The space so  constructed is assigned to
the symbol  {\tt dos3}. The second  homology group of this  space is
verified:

{\tt USER(9): (homology dos3 2) ==>

Computing boundary-matrix in dimension 2.

Rank of the source-module : 1.

Computing boundary-matrix in dimension 3.

Rank of the source-module : 1.

Homology in dimension 2 :

Component Z/2Z

---done---}

Finally we construct the loop  space $\Omega D \Omega S^3$ and compute
the sixth homology group:

{\tt USER(10): (setf odos3 (loop-space dos3)) ==>

[K398 Simplicial-Group]

USER(11): (homology odos3 6) ==>

Computing boundary-matrix in dimension 6.

Rank of the source-module : 26.

Computing boundary-matrix in dimension 7.

Rank of the source-module : 50.

Homology in dimension 6 :

Component Z/6Z

Component Z/2Z

Component Z/2Z

Component Z/2Z

Component Z/2Z

Component Z/2Z

Component Z/2Z

Component Z/2Z

Component Z/2Z

Component Z/2Z

Component Z/2Z

Component Z/2Z

Component Z/2Z

---done---}

\noindent $H_6(\Omega D \Omega S^3) =Z_2^{12} \oplus Z_6$.
The Kenzo program  has determined these groups up to  dimension $9$ in a
few days.

 \vskip .5cm

\centerline{REFERENCES}
\vskip 6pt
1. J.F.Adams. {\it On the cobar construction.} Proc. Nat. Acad. Sci.
42(1956), 409--412.

2. H.J.Baues. {\it The double bar and cobar constructions.} Comp. Math.
43(1981), 331--341.

3. J.P.May. {\it The Geometry of Iterated Loop Spaces.} Lect. Notes in
Math. 1972 v. 271.

4. F.R.Cohen, T.J.Lada, J.P.May. {\it The Homology of Iterated Loop
Spaces.} Lect. Notes in Math. 1976, v. 533.

5. E.Dyer, R.Lashof. {\it Homology of iterated loop spaces.} Amer.
J. Math. 84(1962), 35--88.

6. R.J.Milgram. {\it Iterated loop spaces.} Ann. of Math. 84(1966), N 3,
386--403.

7. V.A.Smirnov. {\it On the cochain complex of topological spaces.}
Mat. Sb. (Russia), 115(1981), 146--158.

8. V.A.Smirnov. {\it Homotopy theory of coalgebras.} Izv. Ac. Nauk
(Russia), 49(1985), 1302--1321.

9. V.A.Smirnov. {\it On the chain complex of an iterated loop
space.} Izv. Ac. Nauk (Russia), 53(1989), 1108--1119.

10. F.Cohen, R.Levi. {\it On the homotopy type of infinite stunted
projective spaces.} To appear in BCAT 1998.

11. F.Sergeraert. {\it The computability Problem in Algebraic
Topology.} Advances in Math. 104(1994), N 1, 1--29.

12. A.Vinogradov, M.Vinogradov.  {\it On  multiple  generalizations of {L}ie
algebras and {P}oisson manifolds.} Secondary calculus and cohomological
physics (Moscow,  1997).  Amer.  Math.  Soc.  Providence, RI, 1998, p.
273--287.

13. J.M.Boardman and R.M.Vogt. {\it Homotopy invariant algebraic structures
on topological spaces.} Lect. Notes in Math. 347, 1973.

14. S.Araki,  T.Kudo.  {\it Topology  of  $H_n$-spaces   and   $H$-squaring
operations.} Mem. Fac. Sci. Kyusyu Univ., Ser. A, 10(1956), N 2, 85--120.

15. F.P.Peterson. {\it Functional cohomology operations.} Trans. Amer. Math.
Soc. 86(1957), p.187--197.

16. V.A.Smirnov.  {\it Functional  homology  operations  and  weak homotopy
type.} Mat. Zametki (Russia), 45(1989), N 5, p.76--86.

17. R. Brown. {\it The twisted Eilenberg-Zilber theorem.}
Celebrazioni Arch. Secolo XX, Simp. Top., 1967,  34--37.

18. G. Carlsson, R. J. Milgram.
{\it Stable homotopy and iterated loop spaces.} In [20],  505--583.

19. V.K.A.M. Gugenheim. {\it On a perturbation theory for the homology of the
loop space.} Journal of Pure and Applied Algebra, 25(1982), 197--205.

20. Handbook of Algebraic Topology (Edited by I.M. James).
North-Holland, 1995.

21. D. M. Kan. {\it A combinatorial definition of homotopy groups.}
Annals of Mathematics.  67(1958),  282--312.

22. J. Rubio, F. Sergeraert.
{\it A program computing the homology groups of loop spaces.}
SIGSAM Bulletin,  25(1991), 20--24.

23. J. Rubio, F. Sergeraert.
{\it Constructive Algebraic Topology.} In preparation, see ftp://
www-fourier.ujf-grenoble.fr/~sergerar.

24. J. Rubio, F. Sergeraert, Y. Siret.
{\it The EAT program.} ftp://www-fourier.ujf-grenoble.fr/\ { }ftp/EAT.

25. F. Sergeraert. {\it The Kenzo program.}
{http://www-fourier.ujf-grenoble. fr/\ { }sergerar/Kenzo/}

\footnote""{
Institut Fourier, BP 74, 38402 St Martin d'H\`eres Cedex, France}

\footnote""{E-mail address: Sergeraert\@ujf-grenoble.fr}

\enddocument